\documentclass[12pt,oneside,reqno,a4paper]{amsart}
\usepackage[utf8]{inputenc}
\usepackage[english]{babel}
\usepackage[T1]{fontenc}
\usepackage{lmodern}
\usepackage{amssymb, amsmath, amsthm, mathtools, cases,enumerate}
\usepackage{color}
\usepackage{mathscinet,geometry}
\usepackage{natbib}
\usepackage{colonequals}
\usepackage{enumitem}
\usepackage[colorlinks=true,linkcolor=blue,urlcolor=blue,citecolor=blue]{hyperref}

\theoremstyle{plain}
\newtheorem{theorem}{Theorem}[section]
\newtheorem{corollary}[theorem]{Corollary}
\newtheorem{proposition}[theorem]{Proposition}

\theoremstyle{definition}
\newtheorem{remark}[theorem]{Remark}
\newtheorem{definition}[theorem]{Definition}

\newcommand{\N}{\mathbb N}

\newcommand{\R}{R}

\DeclareMathOperator*{\esup}{ess\,sup}

\newcommand{\dx}{{\fam0 d}}
\renewcommand{\d}[1]{\,\dx #1}

\newcommand{\M}{\mathfrak{M}}

\newcommand{\Mpl}{\M^+}

\newcommand{\kf}[1][f]{\kappa_{#1}}

\newcommand{\ac}{\mathbin{\ll}}

\begin{document}

\title[RI hulls of weighted Lebesgue spaces]{\textsc{Rearrangement-invariant hulls of weighted Lebesgue spaces}}

\renewcommand{\urladdrname}{\itshape ORCID}
\author[M.~K\v repela]{Martin K\v repela$^*$}
\address{Czech Technical University in Prague, Faculty of Electrical Engineering, Department of Mathematics, Technick\'a~2, 166~27 Praha~6, Czech Republic}
\email{martin.krepela@fel.cvut.cz}
\urladdr{\href{https://orcid.org/0000-0003-0234-1645}{0000-0003-0234-1645}}

\author[Z.~Mihula]{Zden\v ek Mihula$^*$}
\address{Czech Technical University in Prague, Faculty of Electrical Engineering, Department of Mathematics, Technick\'a~2, 166~27 Praha~6, Czech Republic}
\email{mihulzde@fel.cvut.cz}
\urladdr{\href{https://orcid.org/0000-0001-6962-7635}{0000-0001-6962-7635}}

\author[J.~Soria]{Javier Soria$^{\ast\ast}$}
\address{Interdisciplinary Mathematics Institute (IMI), Department of Analysis and Applied Mathematics, Complutense University  of Madrid, 28040 Madrid, Spain}
\email{javier.soria@ucm.es}
\urladdr{\href{https://orcid.org/0000-0003-3098-7056}{0000-0003-3098-7056}}

\thanks{$^*$M.~K\v repela and Z.~Mihula were supported by the project GA23-04720S of the Czech Science Foundation}

\thanks{$^{\ast\ast}$J.~Soria  was partially supported by grants PID2020-113048GB-I00 funded by MCIN/AEI/ 10.13039/501100011033, and Grupo UCM-970966}



\setcitestyle{numbers}
\bibliographystyle{plainnat}

\begin{abstract}
We characterize the rearrangement-invariant hull, with respect to a given measure $\mu$, of weighted Lebesgue spaces. The solution leads us to first consider when this space is contained in the sum of $(L^1 + L^\infty)(\R, \mu)$ and the final condition is given in terms of embeddings for weighted Lorentz spaces.
\end{abstract}
\maketitle

\section{Introduction}
We are interested in studying  the connections between two a priori different theories involving weighted inequalities, namely those related to the Lebesgue spaces   $L^p(\mathbb R^n,u)$ and the weighted Lorentz spaces $\Lambda^p(\mathbb R^n,v)$, intimately connected to the weighted Hardy's inequalities for nonincreasing functions,  and described by the \mbox{$B_p$-class} of Ari\~no and Muckenhoupt \cite{MR989570}.
\medskip

In particular,  we are going to find the rearrangement-invariant (r.i.) hull (the smallest among such spaces) containing $L^p(\mathbb R^n,u)$. To this end, we will actually work on more general measure spaces $(R,\mu)$. 
\medskip

Some related results, dealing with minimal r.i.\  embeddings for other kind  of spaces, can be found in \cite{MR615621}, where the author shows that $L^p(\mathbb R^n)$ is the r.i.\ hull, with respect to the Lebesgue measure, of the Hardy space $H^p(\mathbb R^n)$, $0<p\le 1$. Similarly, in \cite{EKP:00, KP:06} the authors give a characterization of the smallest r.i.\ space for which Sobolev's inequality holds (see also \cite{CPS:15} and references therein). For some other related topics on these questions, see also \cite{MR2262450}.

\medskip

In Section~\ref{prelim}, we will start by providing the tools that we need concerning the nonincreasing and nondecreasing rearrangement of a function. Our main results are Theorem~\ref{thm:our_Ryff} and Proposition~\ref{prop:reverse_H-L}. In Section~\ref{inclusion} we characterize, in Theorem~\ref{embedding}, the embedding of $L^p(\R, \nu)$ in the sum $(L^1 + L^\infty)(\R, \mu)$, which is a necessary condition for the existence of the r.i.\  hull of $L^p(\R, \nu)$, with respect to a different measure $\mu$. In Section~\ref{hullri} we finally describe, in Theorem~\ref{embleblor}, the r.i.\  hull of $L^p(\R, \nu)$. 

\section{Preliminaries}\label{prelim}
Throughout the paper, unless specified otherwise, $(\R, \mu)$ is a $\sigma$-finite nonatomic measure space. We recall the definitions of the \emph{distribution function} of $f$:

\begin{equation*}
\mu_f(s) = \mu(\{x\in \R\colon |f(x)| > s\}),\ s\in(0,\infty),
\end{equation*}
as well as its \emph{nonincreasing rearrangement}
\begin{equation*}
f^*(t) = \inf\{s > 0\colon \mu_f(s) \leq t\},\ t\in(0,\infty).
\end{equation*}
Observe that both $f$ and $f^*$ are equimeasurable; i.e.,  $\mu_f=\mu_{f^*}$ (\cite[Chapter~2, Proposition~1.7]{BS}). 

Let $(\R_1, \mu_1)$ and  $(\R_2, \mu_2)$ be two $\sigma$-finite measure spaces. We say that a function $\sigma\colon \R_1 \to \R_2$ is a \emph{measure-preserving transformation} if, for every $\mu_2$-measurable $E\subset \R_2$, $\sigma^{-1} E = \{x\in \R_1\colon \sigma(x) \in E\}$ is $\mu_1$-measurable and $\mu_1(\sigma^{-1}E) = \mu_2(E)$.

By $\mathfrak M(R,\mu)$, we denote the set of all $\mu$-measurable functions on $R$. The symbol $\mathfrak M^+(R,\mu)$ denotes the set of all nonnegative functions from $\mathfrak M(R,\mu)$.

Let $p \in (0, \infty)$ and  $v\colon(0,\infty) \to [0, \infty)$ be a (Lebesgue-)measurable function. Then  $\Lambda^p(\R, v)$ denotes the set
\begin{equation*}
\Lambda^p(\R, v) = \{f\in\M(\R, \mu)\colon \|f\|_{\Lambda^p(\R, v)} < \infty\},
\end{equation*}
where
\begin{equation*}
\|f\|_{\Lambda^p(\R, v)} = \left(\int_0^\infty f^*(t)^p v(t) \d{t}\right)^{1/p}.
\end{equation*}
Here the rearrangement is taken with respect to $\mu$ and we usually refer to $ \Lambda^p(\R, v)$ as a weighted Lorentz ``space'' \cite{MR2308059}.

Finally, when $A\subset\mathbb R^n$ is a (Lebesgue) measurable set, we denote its Lebesgue measure by $|A|$.

\begin{definition}
Let $f\in\M(\R,\mu)$. We define its \emph{lower distribution function} $\kf\colon (0, \infty) \to [0, \infty]$ as
\begin{equation*}
\kf(s) = \mu(\{x\in \R\colon |f(x)| < s\}),\ s\in(0, \infty).
\end{equation*}

The \emph{nondecreasing rearrangement} $f_*\colon (0, \infty) \to [0, \infty]$ of $f$ is defined as
\begin{equation*}
f_*(t) = \sup\{s > 0\colon \kf(s) < t\},\ t\in(0,\infty).
\end{equation*}
\end{definition}

\begin{remark}\label{rem:another_expression_for_nondecreasing_rearrangement}
Note that, for every $t\in(0,\infty)$,
\begin{equation*}
f_*(t) = \inf\{s > 0\colon \kf(s) \geq t\}.
\end{equation*}
\end{remark}

Some basic properties of the nondecreasing rearrangement, which we shall use later, are listed in the following proposition.
\begin{proposition}
Let $f, g, f_n\in\M(\R, \mu)$, $n\in\N$.
\begin{enumerate}[label=\textnormal{(\roman*)}, leftmargin=25pt]
	\item The lower distribution function $\kf$ is nondecreasing and left-continuous.
	
	\item The nondecreasing rearrangement $f_*$ is nondecreasing and
	\begin{equation}\label{prop:properties_of_inc_rearr:inc_rearr_is_lower_distr}
		f_*(t) = |\{s > 0\colon \kf(s) < t\}| \quad \text{for every $t\in(0, \infty)$}.
	\end{equation}
	In particular, $f_*$ is left-continuous.
	
	\item If $|f| \leq |g|$ $\mu$-a.e., then $\kf\geq \kf[g]$ and $f_* \leq g_*$.
	
	\item If $E\subset \R$ is $\mu$-measurable, then $f_* \leq (f\restriction_E)_*$.
	
	\item If $\mu(\{x\in\R\colon f(x) = 0\}) = \infty$, then $f_*\equiv0$. If $\mu(\R) < \infty$, then $f_*(t) = \infty$ for every $t\in (\mu(\R), \infty)$.
	
	\item If $f_n \downarrow f\geq0$ $\mu$-a.e., then $\kf[f_n]\uparrow \kf$ on $(0, \infty)$.
	
	\item If $f_n \downarrow f\geq0$ $\mu$-a.e., then $(f_n)_* \downarrow f_*$ a.e.~on $(0, \infty)$.
\end{enumerate}
\end{proposition}
\begin{proof}
The monotonicity of both $\kf$ and $f_*$ is obvious. The left-continuity of $\kf$ is an immediate consequence of the monotone convergence theorem. Indeed, if $s_n \uparrow s$, then
\begin{equation*}
\kf(s_n) = \int_{\R} \chi_{\{x\in \R\colon |f(x)| < s_n\}} \d{\mu} \uparrow \int_{\R} \chi_{\{x\in \R\colon |f(x)| < s\}} \d{\mu} = \kf(s).
\end{equation*}
To prove \eqref{prop:properties_of_inc_rearr:inc_rearr_is_lower_distr}, observe that
\begin{equation*}
f_*(t) = \sup\{s \geq 0\colon \kf(s) < t\} \leq |\{s > 0\colon \kf(s) < t\}|
\end{equation*}
inasmuch as $\kf(s) < t$ implies $|\{\tau > 0 \colon \kf(\tau) < t\}|\geq s$. On the other hand, since $\kf(s) \geq t$ implies $|\{\tau > 0 \colon \kf(\tau) < t\}| < s$, we also have
\begin{equation*}
f_*(t) = \inf\{s \geq 0\colon \kf(s) \geq t\} \geq |\{s > 0\colon \kf(s) < t\}|.
\end{equation*}
Hence \eqref{prop:properties_of_inc_rearr:inc_rearr_is_lower_distr} is true. In particular, since $f_*$ is a lower distribution function, it is left-continuous.

The statements (iii), (iv) and (v) are obvious.

Assume that $f_n \downarrow f\geq0$ $\mu$-a.e. Set $F = \{x\in \R\colon \lim_{n\to \infty}f_n(x) = f(x)\}$. Let $s > 0$. Set $E = \{x\in F\colon f(x) < s\}$ and $E_n = \{x\in F\colon f_n(x) < s\}$, $n\in \N$. Clearly $E \subset \bigcup_{n = 1}^\infty E_n$ and $E_1\subset E_2\subset\cdots\subset E$, and so
\begin{equation*}
\mu(E) \leq \mu(\bigcup_{n = 1}^\infty E_n) = \lim_{n \to \infty} \mu(E_n) \leq \mu(E).
\end{equation*}
Since $\mu(\R \setminus F) = 0$, it follows that $\kf(s) = \lim_{n\to\infty} \kf[f_n](s)$. Moreover, the sequence $\{\kf[f_n]\}_{n = 1}^\infty$ is clearly nondecreasing. Hence we have proved (vi). 

It remains to prove (vii); that is, $(f_n)_* \downarrow f_*$ a.e.~on $(0, \infty)$.  The sequence $\{(f_n)_*\}_{n = 1}^\infty$ is clearly nonincreasing. Let $t$ be a point of continuity of $f_*$. We may assume that $f_*(t) < \infty$, for otherwise $(f_n)_*(t) = \infty$, for every $n\in\N$. Clearly
\begin{equation*}
f_*(\tau) - f_*(t) \geq |\{s > 0\colon \kf(s) = t\}|, \quad \text{for every $\tau > t$}.
\end{equation*}
Since $f_*$ is continuous at $t$, it follows that $|\{s > 0\colon \kf(s) = t\}| = 0$. Hence
\begin{equation}\label{prop:properties_of_inc_rearr:inc_rearr_is_lower_distr:eq1}
|\{s > 0\colon \kf(s) < t\}| = |\{s > 0\colon \kf(s) \leq t\}|.
\end{equation}
By $(vi)$, we have $\kf[f_n]\uparrow \kf$ on $(0,\infty)$, and so
\begin{equation*}
\{s > 0\colon \kf(s) > t\} \subset \bigcup_{n = 1}^\infty \{s > 0\colon \kf[f_n](s) > t\}.
\end{equation*}
Consequently,
\begin{align}
\bigcap_{n = 1}^\infty \{s > 0\colon \kf[f_n](s) < t\} &\subset \bigcap_{n = 1}^\infty \{s > 0\colon \kf[f_n](s) \leq t\} \notag\\
&\subset \{s > 0\colon \kf(s) \leq t\}. \label{prop:properties_of_inc_rearr:inc_rearr_is_lower_distr:eq2}
\end{align}
Furthermore,
\begin{equation}\label{prop:properties_of_inc_rearr:inc_rearr_is_lower_distr:eq3}
\{s > 0\colon \kf(s) < t\} \subset \{s > 0\colon \kf[f_n](s) < t\}, \quad \text{for every $n\in\N$},
\end{equation}
since $\kf[f_n]\leq \kf$ for every $n\in\N$. Note that
\begin{equation*}
\Big| \bigcap_{n = 1}^\infty \{s > 0\colon \kf[f_n](s) < t\} \Big| = \lim_{n\to\infty} |\{s > 0\colon \kf[f_n](s) < t\}|,
\end{equation*}
even when $|\{s > 0\colon \kf[f_n](s) < t\}| = \infty$ for every $n\in\N$. Indeed, this follows from the observation that $|\{s > 0\colon \kf[f_n](s) < t\}| = \infty$ if and only if $\{s > 0\colon \kf[f_n](s) < t\} = (0, \infty)$. Therefore, using \eqref{prop:properties_of_inc_rearr:inc_rearr_is_lower_distr:eq1}, \eqref{prop:properties_of_inc_rearr:inc_rearr_is_lower_distr:eq2} and \eqref{prop:properties_of_inc_rearr:inc_rearr_is_lower_distr:eq3}, we arrive at
\begin{align*}
f_*(t) &= |\{s > 0\colon \kf(s) < t\}| = |\{s > 0\colon \kf(s) \leq t\}| \\
&\geq \Big|\bigcap_{n = 1}^\infty \{s > 0\colon \kf[f_n](s) < t\}\Big| = \lim_{n\to\infty} |\{s > 0\colon \kf[f_n](s) < t\}| \\
&\geq |\{s > 0\colon \kf(s) < t\}| = f_*(t).
\end{align*}
Hence $f_*(t) = \lim_{n\to\infty} |\{s > 0\colon \kf[f_n](s) < t\}| = \lim_{n\to\infty} (f_n)_*(t)$.
\end{proof}

\begin{remark}\label{rem:nondecreasing_by_means_of_nonincreasing}
Assume that $\mu(\R) < \infty$, and let $f\in\M(\R,\mu)$. Note that
\begin{align*}
f^*(\mu(\R) - t) &= \inf\{s > 0\colon \mu_f(s) \leq \mu(\R) - t\} \\
&= \inf\{s > 0\colon \mu(\{x\in \R\colon |f(x)| \leq s\}) \geq t\} \\ 
&= \inf\{s > 0\colon \kf(s) \geq t\} = f_*(t),
\end{align*}
for every $t\in(0,\mu(\R))$. Therefore, our definition of the nondecreasing rearrangement coincides with its usual definition on finite measure spaces (e.g., see \cite[p.~78]{B:19}).
\end{remark}

The following result is a natural extension of Ryff's theorem \cite{R:70} to the case of the nondecreasing rearrangement (see also \cite{CafSor}):

\begin{theorem}\label{thm:our_Ryff}
Let $f\in\Mpl(\R, \mu)$. Set  $T = \esup_{x\in\R} f(x)$. Then there exists a measure-preserving transformation $\sigma\colon \R \to (0, \mu(\R))$ such that
\begin{equation*}
f = f_* \circ \sigma, \quad \text{$\mu$-a.e.\ in $\R$},
\end{equation*}
if and only if one of the following conditions is satisfied:
\begin{enumerate}[label=\textnormal{(\roman*)}, leftmargin=20pt]
	\item $\kf(T) < \infty$, 
	
	\item $\kf(s) < \infty$, for every $s\in(0, T)$, and $\mu(\{x\in\R\colon f(x) = T\}) = 0$.
\end{enumerate}
\end{theorem}
\begin{proof}
First, assume that $\mu(\R) < \infty$. Note that the condition (i) is plainly always satisfied. By the standard Ryff theorem \cite{R:70}, there is a measure-preserving transformation $\sigma_1\colon \R \to (0, \mu(R))$ such that $f = f^* \circ \sigma_2$, $\mu$-a.e. Owing to Remark~\ref{rem:nondecreasing_by_means_of_nonincreasing}, we have $f^* = f_* \circ \sigma_2$, where $\sigma_2\colon (0, \mu(\R)) \to (0, \mu(\R))$ is defined as $\sigma_2(t) = \mu(\R) - t$, $t\in (0,\mu(\R))$. It is easy to see that $\sigma = \sigma_2 \circ \sigma_1$ is the desired measure-preserving transformation.

Second, assume that $\mu(\R) = \infty$ and that one of the conditions holds. Redefining $f$ on a set of measure $0$ if necessary, we may assume that $\{x\in\R\colon f(x) > T\} = \emptyset$. For $n\in\N$, we define numbers $T_n$ as
\begin{equation*}
T_n = \begin{cases}
		
		T(1 - \frac1{n}), \quad&\text{if $T< \infty$},\\
		
		n-1, \quad&\text{if $T = \infty$}.
		
\end{cases}
\end{equation*}
We also set
\begin{align*}
\R_n &= \{x\in\R\colon T_n\leq f(x) < T_{n+1}\},\ n\in\N, \\
\R_\infty &= \{x\in\R\colon f(x) = T\}.
\end{align*}
Clearly
\begin{equation*}
\R = \bigcup_{n = 1}^\infty R_n \cup \R_\infty,
\end{equation*}
and the sets $\R_n$, $n\in\N$, and $\R_\infty$ and mutually disjoint. We have
\begin{equation*}
\mu(\R_n) = \kf(T_{n+1}) - \kf(T_n) < \infty, \quad \text{for every $n\in\N$}.
\end{equation*}
Set
\begin{equation*}
f_n = f\restriction_{\R_n},\ n\in\N.
\end{equation*}
Thanks to the first part of this proof, there are measure-preserving transformations $\widetilde{\sigma}_n\colon \R_n \to (0, \mu(\R_n))$ such that
\begin{equation*}
f_n = (f_n)_*\circ \widetilde{\sigma}_n, \quad \text{$\mu$-a.e.~in $\R_n$}.
\end{equation*}
Set
\begin{equation*}
\sigma_n(x) = \kf(T_n) + \widetilde{\sigma}(x),\ x\in\R_n.
\end{equation*}
Clearly, $\sigma_n\colon \R_n \to (\kf(T_n), \kf(T_n) + \mu(\R_n))$ is measure-preserving. Note that the intervals $\{(\kf(T_n), \kf(T_n) + \mu(\R_n))\}_{n=1}^\infty$ are disjoint, $\kf(T_1) = 0$, and $\kf(T_n) = \kf(T_{n-1}) + \mu(\R_{n-1})$, $n\geq 2$. Furthermore, we clearly have $\lim_{n\to\infty}\kf(T_n) = \kf(T)$. Therefore
\begin{equation*}
(0, \infty) = \bigcup_{n = 1}^\infty I_n \cup (\kf(T), \infty) \cup S,
\end{equation*}
where $I_n$ are the intervals $(\kf(T_n), \kf(T_n) + \mu(\R_n))$ and $S$ is a countable set. Moreover, the sets on the right-hand side are mutually disjoint. Since
\begin{equation*}
\kf(s) = \kf[f_n](s) + \kf(\min\{s, T_n\}) + \mu(\{x\in\R\colon T_{n+1}\leq f(x) < s\}),
\end{equation*}
for every $n\in\N$ and $s>0$, it follows that
\begin{equation*}
f(x) = f_n(x) = f_*(\widetilde{\sigma}_n(x) + \kf(T_n)) = f_*(\sigma_n(x)), \quad \text{$\mu$-a.e.~in $R_n$}.
\end{equation*}
Now, assume that (ii) holds. Redefining $f$ on a set of measure $0$ if necessary, we may assume that $\R_\infty = \emptyset$. Note that $\kf(T) = \infty$, for otherwise $\mu(R) = \kf(T) < \infty$. Set
\begin{equation*}
\sigma(x) = \sum_{n = 1}^\infty \sigma_n(x)\chi_{\R_n}(x),\ x\in\R.
\end{equation*}
Clearly, for $\mu$-a.e.~every $x\in\R$,
\begin{equation*}
f(x) = f_n(x) = f_*(\sigma_n(x)) = f_*(\sigma(x)).
\end{equation*}
Here $n$ is the unique $n\in\N$ such that $x\in\R_n$. Hence $f = f_*\circ \sigma$ $\mu$-a.e.~in $\R$. We need to show that $\sigma\colon \R \to (0, \infty)$ is a measure-preserving transformation. Let $E\subset (0, \infty)$ be measurable. We have
\begin{align*}
\sigma^{-1}E &= \bigcup_{n = 1}^\infty \sigma^{-1}(E\cap I_n) \cup  \sigma^{-1}(E\cap S) = \bigcup_{n = 1}^\infty \sigma_n^{-1}(E\cap I_n) \cup \emptyset \\
& = \bigcup_{n = 1}^\infty \sigma_n^{-1}(E\cap I_n).
\end{align*}
Since the sets $\sigma_n^{-1}(E\cap I_n)$ are $\mu$-measurable inasmuch as the functions $\sigma_n\colon \R_n \to I_n$ are measure-preserving, $\sigma^{-1}E$ is $\mu$-measurable. Furthermore, since $(0, \infty) = \bigcup_{n = 1}^\infty I_n$ up to a countable set, the intervals $I_n$ are mutually disjoint, and so are the sets $R_n$, we have
\begin{equation*}
\mu(\sigma^{-1}E) = \sum_{n = 1}^\infty \mu(\sigma_n^{-1}(E\cap I_n)) = \sum_{n = 1}^\infty |E\cap I_n| = |E|.
\end{equation*}
Here we used the fact that $\sigma_n\colon \R_n \to I_n$ are measure-preserving again. Hence $\sigma$ is the desired measure-preserving transformation if (ii) holds. 

Now assume that (i) holds. The proof in this case follows along the same lines as that in the previous case, and so we just sketch it. One has $\R_\infty = \infty$, for otherwise $\mu(\R) = \kf(T) + \mu(\R_\infty) < \infty$. Since $\mu$ is nonatomic, $\sigma$-finite and $\mu(\R_\infty) = \infty$, it follows from \cite[Chapter~2, Proposition~7.4]{BS} that there is a measure-preserving transformation $\sigma_\infty\colon \R_\infty \to (\kf(T), \infty)$. Set
\begin{equation*}
\sigma(x) = \sum_{n = 1}^\infty \sigma_n(x)\chi_{\R_n}(x) + \sigma_\infty(x)\chi_{\R_\infty}(x),\ x\in\R.
\end{equation*}
As in the case (ii), $\sigma\colon \R \to (0, \infty)$ is a measure-preserving transformation. Indeed, we have
\begin{equation*}
\sigma^{-1}E = \bigcup_{n = 1}^\infty \sigma_n^{-1}(E\cap I_n) \cup \sigma_\infty^{-1}(E\cap (\kf(T), \infty)),
\end{equation*}
for every measurable $E\subset (0, \infty)$, and
\begin{align*}
\mu(\sigma^{-1}E) &= \sum_{n = 1}^\infty \mu(\sigma_n^{-1}(E\cap I_n)) + \mu(\sigma_\infty^{-1}(E\cap (\kf(T), \infty))) \\
&= \sum_{n = 1}^\infty |E\cap I_n| + |E\cap (\kf(T), \infty)| = |E|.
\end{align*}
Finally, note that $f_*(t) = T$, for every $t\in (\kf(T), \infty)$. Therefore, if $x\in\R_\infty$, then
\begin{equation*}
f(x) = T = f_*(\sigma_\infty(x)).
\end{equation*}
As in the case (ii), we have
\begin{equation*}
f(x) = f_*(\sigma(x)), \quad \text{for $\mu$-a.e.~$x\in\R\setminus\R_\infty$}.
\end{equation*}
Hence $f = f_*\circ \sigma$ $\mu$-a.e.~in $\R$, and so $\sigma$ is the desired measure-preserving transformation if (i) holds.

Last, we turn our attention to the necessity of the conditions (i) and (ii). Note that, if neither of them is satisfied, then there is $s\in(0, T]$ such that
\begin{equation*}
\kf(s) = \infty \quad \text{and} \quad \mu(\{x\in\R\colon f(x) \geq s\}) > 0.
\end{equation*}
Here we used the fact that $T$ is the essential supremum of $f$ over $\R$. Since $\infty = \kf(s) = \lim_{n\to \infty} \kf(s_n)$ for every $0 <s_n \uparrow s$, we have
\begin{equation*}
f_*(t) < s, \quad \text{for all $t\in(0, \infty)$}.
\end{equation*}
Hence, for every $x\in \{y\in \R\colon f(y) \geq s\}$ and any function $\sigma\colon \R \to (0, \infty)$, 
\begin{equation*}
f(x) \geq s > f_*(\sigma(x)).
\end{equation*}
Therefore, it is not possible to find a (measure-preserving transformation) $\sigma\colon \R \to (0, \infty)$ such that $f = f_* \circ \sigma$ $\mu$-a.e.~in $\R$, inasmuch as $\mu(\{x\in \R\colon f(x) \geq s\}) > 0$.
\end{proof}

The standard Hardy--Littlewood inequality for nonincreasing rearrangements tells us that, for every $f,g\in\Mpl(\R, \mu)$,
\begin{equation*}
\int_{\R} f g \d{\mu} \leq \int_0^{\mu(\R)} f^*(t) g^*(t) \d{t}.
\end{equation*}
As a corollary, when $\mu(\R) < \infty$, it is easy to obtain (e.g., \cite[Corollary~2.18]{B:19}) its reverse version; namely
\begin{equation*}
\int_{\R} f g \d{\mu} \geq \int_0^{\mu(\R)} f^*(t) g_*(t) \d{t}.
\end{equation*}
The following proposition extends this reverse Hardy--Littlewood inequality to infinite measure spaces.
\begin{proposition}\label{prop:reverse_H-L}
Let $f,g\in\Mpl(\R, \mu)$. We have
\begin{equation*}
\int_{\R} f g \d{\mu} \geq \int_0^{\mu(\R)} f^*(t) g_*(t) \d{t}.
\end{equation*}
\end{proposition}
\begin{proof}
We start with the case $f\equiv1$; that is, we first show that
\begin{equation}\label{prop:reverse_H-L:f==1}
\int_{\R} g \d{\mu} \geq \int_0^{\mu(\R)} g_*(t) \d{t}.
\end{equation}
Let $\R_1 \subset \R_2 \subset \cdots \subset \R$ be $\mu$-measurable sets such that $\R = \bigcup_{n = 1}^\infty \R_n$ and \mbox{$\mu(\R_n) < \infty$}, for every $n\in\N$. Set
\begin{equation*}
g_n = g\restriction_{\R_n},\ n\in\N.
\end{equation*}
For every $n\in\N$, we define functions $h_n\colon (0, \infty) \to [0, \infty]$ as
\begin{equation*}
h_n(t) = \inf\{s>0 \colon \kf[g_n](s) \geq t\},\ t\in(0, \infty).
\end{equation*}
Clearly $g_*(t) \leq h_n(t)$ for every $t\in (0, \infty)$ and $n\in\N$. Furthermore, by Remark~\ref{rem:another_expression_for_nondecreasing_rearrangement} and \ref{rem:nondecreasing_by_means_of_nonincreasing}, we have
\begin{equation*}
h_n(t) = (g_n)^*(\mu(\R_n) - t), \quad \text{for every $t\in(0,\mu(\R_n))$}.
\end{equation*}
Therefore, using the layer cake representation (e.g., see \cite[Chapter~2, Proposition~1.8]{BS}), we have
\begin{align*}
\int_0^{\mu(\R_n)} g_*(t) \d{t} &\leq \int_0^{\mu(\R_n)} h_n(t) \d{t} = \int_0^{\mu(\R_n)} (g_n)^*(\mu(\R_n) - t) \d{t} \\
&= \int_0^{\mu(\R_n)} (g_n)^*(t) \d{t} = \int_{\R_n} g_n \d{\mu} = \int_{\R_n} g \d{\mu}.
\end{align*}
Hence, letting $n\to\infty$, we obtain \eqref{prop:reverse_H-L:f==1}.

Next, let $f\in\Mpl(\R, \mu)$ be a simple function having the form
\begin{equation*}
f = \sum_{i = 1}^n \alpha_i \chi_{E_i},
\end{equation*}
where $\alpha_i\geq0$ and $E_i$ are $\mu$-measurable sets such that $E_1\subset E_2 \subset \cdots \subset E_n$ and $\mu(E_i) < \infty$. Note that
\begin{equation*}
f^* = \sum_{i = 1}^n \alpha_i \chi_{(0, \mu(E_i))}.
\end{equation*}
Using \eqref{prop:reverse_H-L:f==1}, we obtain
\begin{align*}
\int_{\R} f g \d{\mu} &= \sum_{i = 1}^\infty \alpha_i \int_{E_i} g \d{\mu} \geq \sum_{i = 1}^\infty \alpha_i \int_0^{\mu(E_i)} (g\restriction_{E_i})_*(t) \d{t} \\
&\geq \sum_{i = 1}^\infty \alpha_i \int_0^{\mu(E_i)} g_*(t) \d{t} = \int_0^{\mu(\R)} \left(\sum_{i = 1}^\infty \alpha_i \chi_{(0, E_i)}(t)\right) g_*(t) \d{t} \\
&= \int_0^{\mu(\R)} f^*(t) g_*(t) \d{t}.
\end{align*}

Finally, for every $f\in\Mpl(\R, \mu)$ there is a sequence $\{f_n\}_{n = 1}^\infty$ of simple functions such that $0\leq f_n \uparrow f$ $\mu$-a.e.~in $\R$. Since we also have $(f_n)^* \uparrow f^*$, the general case follows immediately from the previous one thanks to the monotone convergence theorem.
\end{proof}

\section{Inclusion and endpoint spaces}\label{inclusion}
Let $\mu, \nu$ be two measures on $\R$ over the same $\sigma$-algebra. When $X(\R, \nu)$ and $Y(\R, \mu)$ are Banach spaces of measurable functions in which two functions coinciding $\nu$ and $\mu$-a.e., respectively, are identified, $X(\R, \nu) \hookrightarrow Y(\R,\mu)$ denotes the following fact. The operator $I\colon X(\R, \nu) \to Y(\R,\mu)$ defined as
\begin{equation*}
I([f]_{\nu\text{-a.e.}}) = [f]_{\mu\text{-a.e.}},f \in X(\R, \nu),
\end{equation*}
is well defined and bounded. We will denote as $\|X\hookrightarrow Y\|$ the norm of the inclusion operator $I$. Note that, should $I$ be well defined, every two functions that coincide $\nu$-a.e.\ have to coincide $\mu$-a.e., too. In other words, $\mu$ has to be absolutely continuous with respect to $\nu$ (cf.~\cite[p.~342]{K:81}). We shall denote this fact by $\mu \ac \nu$, as usual. Finally, note that the operator $I$ need not be injective unless we also have $\nu \ac \mu$.

\begin{theorem}\label{embedding}
Let $\mu, \nu$ be nonatomic $\sigma$-finite measures on $\R$ over the same \mbox{$\sigma$-algebra}. Let $p\in [1, \infty)$. The following two statements are equivalent.
\begin{enumerate}[label=\textnormal{(\roman*)}, leftmargin=20pt]
\item
\begin{equation}\label{thm:inclusion_to_L1+Linfty:embedding}
L^p(\R, \nu) \hookrightarrow (L^1 + L^\infty)(\R, \mu).
\end{equation}

\item $\mu \ac \nu$ and
\begin{equation}\label{charctnorm}
A = \sup_{\substack{E\subset \R \\ \mu(E) \le 1}} \bigg\|\frac{\dx\mu}{\dx\nu}\bigg\|_{L^{p'}(E, \nu)} < \infty.
\end{equation}
\end{enumerate}

Moreover, if \eqref{thm:inclusion_to_L1+Linfty:embedding} is valid, then we have
\begin{equation*}
\|L^p(\R, \nu) \hookrightarrow (L^1 + L^\infty)(\R, \mu)\| = A.
\end{equation*}
\end{theorem}
\begin{proof}
The rearrangements in this proof are taken with respect to $\mu$.

Assume that \eqref{thm:inclusion_to_L1+Linfty:embedding} is valid. As was pointed out before the theorem, we necessarily have $\mu \ac \nu$. Let $C\in(0, \infty)$ be the optimal constant of the embedding, that is, $C = \|L^p(\R, \nu) \hookrightarrow (L^1 + L^\infty)(\R, \mu)\|$. Let $E\subset \R$ be $\mu$-measurable, with $\mu(E) \leq 1$, and let $f$ be a function such that $f\chi_E\in(L^1 + L^\infty)(\R, \mu)$. We have
\begin{equation*}
\|f\chi_E\|_{L^1(\R, \mu)} = \int_0^{\mu(E)} (f\chi_E)^*(t) \d{t} = \|f\chi_E\|_{(L^1 + L^\infty)(\R, \mu)}.
\end{equation*}
Therefore, using \eqref{thm:inclusion_to_L1+Linfty:embedding}, we obtain
\begin{align*}
\|f\chi_E\|_{L^1(\R, \mu)}    
\leq C  \|f\chi_E\|_{L^p(\R, \nu)}.
\end{align*}
Hence $\|L^p(E, \nu) \hookrightarrow L^1(E, \mu)\|\leq C$. By \cite[Lemmas~3 and 4]{K:81} (cf.~\cite[Theorem~2.1.8]{GreyPhD}), we have
\begin{equation*}
\Big\| \frac{\dx\mu}{\dx\nu} \Big\|_{L^{p'}(E, \nu)} \leq C.
\end{equation*}
Thus,  $A \leq C$.
\medskip

Assume now that $\mu \ac \nu$ and $A < \infty$. Let $f\in\M(\R, \mu)$. By \cite[Chapter~2, Proposition~3.3]{BS} and the H\"older inequality, we have
\begin{align*}
\|f\|_{(L^1 + L^\infty)(\R, \mu)} &= \int_0^{\min\{\mu(R), 1\}} f^*(t) \d{t} = \sup_{\mu(E) = \min\{\mu(R), 1\}} \int_E |f| \d{\mu} \\
&= \sup_{\mu(E) = \min\{\mu(R), 1\}} \int_E |f| \frac{\dx\mu}{\dx\nu}\d{\nu} \leq A \|f\|_{L^p(\R,\nu)}.
\end{align*}
Hence $\|L^p(\R, \nu) \hookrightarrow (L^1 + L^\infty)(\R, \mu)\| \leq A$.
\end{proof}

\section{Rearrangement-invariant hulls}\label{hullri}

The following result describes, with full generality, the existence of the minimal r.i. hull containing a Lebesgue space, with a very precise control of the norm estimates. We start by defining such a set:

\begin{definition}
Let  $(\R, \mu)$ be a measure space.

\begin{enumerate}[label=\textnormal{(\roman*)}, leftmargin=20pt]
\item We say that $X\subset \mathfrak{M}(\R, \mu)$ is a rearrangement-invariant (r.i.) set if 
for every $f\in X$ and $g\in \mathfrak{M}(\R, \mu)$ satisfying $g^*=f^*$, one has $g\in X$.
\item Given  $Y\subset  \mathfrak{M}(\R, \mu)$, we say that $X$ is  its   r.i.\  hull, if $Y\subset X$, $X$ is an r.i.\  set and, for every $f\in X$, there exists $g\in Y$ such that $g^*=f^*$.
\end{enumerate}
\end{definition}

\begin{remark}
We observe that if $X$ is the r.i.\  hull of $Y$, and $Z$ is another r.i.\  set with $Y\subset Z$, then $X\subset Z$; i.e., $X$ is the smallest r.i.\  set containing $Y$. In fact, if we pick $f\in X$, then there exists $g\in Y$ such that $g^*=f^*$, which implies that $g\in  Z$, and finally, since $Z$ is  an r.i.\  set, we conclude that $f\in Z$.

\medskip

Another interesting remark is that this definition is really stronger than just assuming minimality in the inclusion. In fact, if $Y=\mathcal{C}[0,1]$, then $L^\infty[0,1]$  is clearly the smallest  r.i.\ space containing  $Y$, but $L^\infty[0,1]$  is not the r.i.\ hull of $Y=\mathcal{C}[0,1]$. Indeed, given $f=\chi_{[0,1/2]}\in L^\infty[0,1]$, there is no $g\in \mathcal{C}[0,1]$ such that $f^*=g^*$.

\end{remark}

In the following, we again assume that $(R,\mu)$ is a~$\sigma$-finite nonatomic measure space.

\begin{theorem}\label{embleblor}
Let $p \in (0 ,\infty)$ and $v\colon \R \to [0, \infty)$ be a $\mu$-measurable function. The following facts are true\textemdash the rearrangements are taken with respect to $\mu$.
\begin{enumerate}[label=\textnormal{(\roman*)}, leftmargin=25pt]
\item We have
\begin{equation}\label{thm:hull_of_weighted_Lp:embedding}
\|f\|_{\Lambda^p(\R, v_*)} \leq \|f\|_{L^p(\R, v)} \quad \text{for every $f\in L^p(\R, v)$}.
\end{equation}

\item If $v_*\not\equiv0$, then, for every $g\in \Lambda^p(\R, v_*)$ and $\varepsilon > 0$, there is a function $f\in L^p(\R, v)$ such that $f^* = g^*$ and
\begin{equation}\label{thm:hull_of_weighted_Lp:almost_equal_norms}
\|g\|_{\Lambda^p(\R, v_*)} \leq \|f\|_{L^p(\R, v)} \leq (1 + \varepsilon) \|g\|_{\Lambda^p(\R, v_*)}.
\end{equation}
Therefore, $\Lambda^p(\R, v_*)$ is the r.i.\ hull of $L^p(\R, v)$. Moreover, if either of the conditions (i), (ii) in Theorem~\ref{thm:our_Ryff} is satisfied with $f$ replaced by $v$, then it is possible to take $\varepsilon = 0$ in \eqref{thm:hull_of_weighted_Lp:almost_equal_norms}. In particular, that is possible when  $\mu(\R) < \infty$.

\item If $v_*\equiv0$, then, for every $g\in\Mpl(\R,\mu)$ such that $g^*(t) < \infty$ for every $t\in(0,\infty)$, and every $\varepsilon > 0$, there is a function $f\in L^p(\R, v)$ such that $f^* = g^*$ and
\begin{equation}\label{thm:hull_of_weighted_Lp:epsilon_norm}
\|f\|_{L^p(\R, v)} \leq \varepsilon.
\end{equation}
\end{enumerate}
\end{theorem}

\begin{proof}
As for (i), the validity of \eqref{thm:hull_of_weighted_Lp:embedding}, follows immediately from Proposition~\ref{prop:reverse_H-L}. In the following, set $T = \esup_{\R} v$.

As next, we are going to prove (ii). Assume that $v_*\not\equiv0$. This implies that the set $\{s\in(0, T]\colon \kf[v](s) < \infty\}$ is nonempty. Set
\begin{align*}
S &= \sup \{s\in(0, T]\colon \kf[v](s) < \infty\},\\
\R_1 &= \{x\in\R\colon v(x) < S\},\\
\intertext{and}
v_1 &= v\restriction_{\R_1}.
\end{align*}
Since $\kf[v](s) < \infty$ for every $s\in(0, S)$ and $\R_1 \cap \{x\in\R\colon v(x) \geq S\} = \emptyset$, there is a measure-preserving transformation $\sigma_1\colon \R_1 \to (0, \mu(\R_1))$ such that
\begin{equation*}
v_1 = (v_1)_* \circ \sigma_1 \quad \text{$\mu$-a.e.\ in $\R_1$}
\end{equation*}
thanks to Theorem~\ref{thm:our_Ryff}. Moreover, since
\begin{equation*}
(v_1)_*(t) = v_*(t) \quad \text{for every $t\in(0, \mu(\R_1))$},
\end{equation*}
we have
\begin{equation}\label{thm:hull_of_weighted_Lp:v=v_down_star_circ_sigma1}
v = v_* \circ \sigma_1 \quad \text{$\mu$-a.e.\ in $\R_1$}.
\end{equation}
We now distinguish between three cases.

First, assume that $\kf[v](S) = \infty$. It follows that $\mu(\R) = \infty$ inasmuch as $\mu(\R_1) = \kf[v](S) = \infty$. Let $g\in \Lambda^p(\R, v_*)$. Set
\begin{equation*}
f(x) = \begin{cases}
	g^*(\sigma_1(x)),\ &x\in\R_1,\\
	0,\ &x\in\R\setminus\R_1.
\end{cases}
\end{equation*}
We claim that the functions $f$ and $g$ are equimeasurable. Indeed, using the fact that the functions $g^*\circ \sigma_1$ and $g^*$ are equimeasurable (\cite[Chapter~2, Proposition~7.2]{BS}) (and so are the functions $g$ and $g^*$), we have
\begin{equation*}
\mu_f(\lambda) = \mu(\{x\in\R_1\colon g^*(\sigma_1(x)) > \lambda \}) = |\{t>0\colon g^*(t) > \lambda\}| = \mu_g(\lambda)
\end{equation*}
for every $\lambda > 0$. Furthermore, we have
\begin{align*}
\int_{\R} f(x)^p v(x) \d{\mu(x)} &= \int_{\R_1} g^*(\sigma_1(x))^p v_*(\sigma_1(x)) \d{\mu(x)}\\
&= \int_0^\infty g^*(t)^p v_*(t) \d{t}.
\end{align*}
Here we used \eqref{thm:hull_of_weighted_Lp:v=v_down_star_circ_sigma1} in the first equality and the fact that $\sigma_1\colon\R_1 \to (0, \mu(\R_1))$ is a measure-preserving transformation in the second one. Hence \eqref{thm:hull_of_weighted_Lp:almost_equal_norms} with $\varepsilon = 0$ follows.

Second, assume that $S < \infty$ and $\kf[v](S) < \infty$. Let $g\in \Lambda^p(\R, v_*)$ and $\varepsilon > 0$. Set
\begin{equation}\label{thm:hull_of_weighted_Lp:R2_def}
\R_2 = \{x\in\R\colon S\leq v(x) < (1 + \varepsilon)^p S\}.
\end{equation}
Note that
\begin{equation}\label{thm:hull_of_weighted_Lp:measure_R=measure_R1+measure_R2}
\mu(\R) = \mu(\R_2) + \kf[v](S).
\end{equation}
Indeed, if $\mu(\R) < \infty$, then $S = T$ and $v(x) = S$ for $\mu$-a.e.~$x\in\R\setminus\R_1$; whence \eqref{thm:hull_of_weighted_Lp:measure_R=measure_R1+measure_R2} follows. If $\mu(\R) = \infty$, then the definition of $S$ and the assumptions $S < \infty$ and $\kf[v](S) < \infty$ imply that $\kf[v]((1 + \varepsilon)^p S) = \infty$; therefore
\begin{equation*}
\mu(\R) = \infty = \infty - \kf[v](S) + \kf[v](S) = \mu(\R_2) + \kf[v](S).
\end{equation*}
Now, it follows from \cite[Chapter~2, Proposition~7.4]{BS} that there is a measure-preserving transformation $\sigma_2\colon \R_2 \to (0, \mu(\R_2))$. Set
\begin{equation*}
f(x) = \begin{cases}
	g^*(\sigma_1(x)),\ &x\in\R_1,\\
	g^*(\kf[v](S) + \sigma_2(x)),\ &x\in\R_2,\\
	0,\ &x\in\R\setminus(\R_1 \cup \R_2).
\end{cases}
\end{equation*}
Note that the functions $f$ and $g$ are equimeasurable. Indeed, we have
\begin{align*}
\mu_f(\lambda) &= \mu(\{x\in\R_1\colon g^*(\sigma_1(x)) > \lambda \}) + \mu(\{x\in\R_2\colon g^*(\kf[v](S) + \sigma_2(x)) > \lambda \}) \\
&= |\{t\in (0, \kf[v](S))\colon g^*(t) > \lambda \}| + |\{t\in (\kf[v](S), \mu(R))\colon g^*(t) > \lambda \}| \\
&= |\{t > 0\colon g^*(t) > \lambda \}| = \mu_g(\lambda)
\end{align*}
for every $\lambda > 0$. Here we used the fact that $\sigma_1$ and $\kf[v](S) + \sigma_2$ are measure-preserving transformations from $\R_1$ and $\R_2$ to $(0, \kf[v](S))$ and $(\kf[v](S), \mu(R))$, respectively. Clearly
\begin{align}
\int_{\R}|f(x)|^p v(x) \d{\mu(x)} &=  \int_{\R_1}|f(x)|^p v(x) \d{\mu(x)} + \int_{\R_2}|f(x)|^p v(x) \d{\mu(x)} \notag\\
&= \int_{\R_1} g^*(\sigma_1(x))^p v_*(\sigma_1(x)) \d{\mu(x)} \notag\\
&\quad + \int_{\R_2} g^*(\kf[v](S) + \sigma_2(x))^p v(x) \d{\mu(x)} \equalscolon I_1 + I_2. \label{thm:hull_of_weighted_Lp:I=I1+I2}
\end{align}
Here we used \eqref{thm:hull_of_weighted_Lp:v=v_down_star_circ_sigma1} in the second equality. Using the fact that $\sigma_1$ is a measure-preserving transformation from $\R_1$ to $(0, \kf[v](S))$, we have
\begin{equation*}
I_1 = \int_0^{\kf[v](S)} g^*(t)^p v_*(t) \d{t}.
\end{equation*}
As for $I_2$, note that the definition of $S$ implies that
\begin{equation*}
v_*(t) = S \quad \text{for every $t\in(\kf[v](S), \mu(\R))$}.
\end{equation*}
On the one hand, using the definition of $\R_2$ and \eqref{thm:hull_of_weighted_Lp:measure_R=measure_R1+measure_R2}, we have
\begin{align}
\int_{\kf[v](S)}^{\mu(\R)} g^*(t)^p v_*(t) \d{t} &= S \int_{\kf[v](S)}^{\mu(\R)} g^*(t)^p \d{t} = S \int_{\R_2} g^*(\kf[v](S) + \sigma_2(x))^p \d{\mu(x)} \notag\\
&\leq \int_{\R_2} g^*(\kf[v](S) + \sigma_2(x))^p v(x) \d{\mu(x)} = I_2. \label{thm:hull_of_weighted_Lp:lower_bound}
\end{align}
On the other hand, we also have
\begin{align}
I_2 &= \int_{\R_2} g^*(\kf[v](S) + \sigma_2(x))^p v(x) \d{\mu(x)} \notag\\
&\leq (1+\varepsilon)^p S \int_{\R_2} g^*(\kf[v](S) + \sigma_2(x))^p \d{\mu(x)} \notag\\
&= (1+\varepsilon)^p S \int_{\kf[v](S)}^{\mu(\R)} g^*(t)^p \d{t} = (1+\varepsilon)^p \int_{\kf[v](S)}^{\mu(\R)} g^*(t)^p v_*(t) \d{t}. \label{thm:hull_of_weighted_Lp:upper_bound}
\end{align}
Here we used the fact that $\kf[v](S) + \sigma_2$ is a measure-preserving transformation from $\R_2$ to $(\kf[v](S), \mu(R))$. Combining \eqref{thm:hull_of_weighted_Lp:I=I1+I2} with \eqref{thm:hull_of_weighted_Lp:lower_bound} and \eqref{thm:hull_of_weighted_Lp:upper_bound}, we arrive at
\begin{align*}
&\int_0^{\kf[v](S)} g^*(t)^p v_*(t) \d{t} + \int_{\kf[v](S)}^{\mu(\R)} g^*(t)^p v_*(t) \d{t} \leq \int_{\R}|f(x)|^p v(x) \d{\mu(x)} \\
&\quad\leq \int_0^{\kf[v](S)} g^*(t)^p v_*(t) \d{t} + (1+\varepsilon)^p \int_{\kf[v](S)}^{\mu(\R)} g^*(t)^p v_*(t) \d{t},
\end{align*}
whence \eqref{thm:hull_of_weighted_Lp:almost_equal_norms} immediately follows.

Third, assume that $S = \infty$ and $\kf[v](S) < \infty$. Note that we necessarily have $S=T$. It follows that $\mu(\R) < \infty$. Indeed, if we had $\mu(\R) = \infty$, our current assumptions would imply that
\begin{equation*}
\mu(\{x\in\R\colon v(x) = \infty\}) = \infty;
\end{equation*}
however, this is not possible. As $\mu(\R) < \infty$, there is a measure-preserving transformation $\sigma\colon\R \to (0, \mu(\R))$ such that
\begin{equation}\label{thm:hull_of_weighted_Lp:v=v_down_star_circ_sigma_finite_measure}
v = v_* \circ \sigma \quad \text{$\mu$-a.e.~in $\R$}
\end{equation}
thanks to Theorem~\ref{thm:our_Ryff}. Let $g\in \Lambda^p(\R, v_*)$. Set
\begin{equation*}
f(x) = g^*(\sigma(x)),\ x\in\R.
\end{equation*}
Using the fact that the functions $g^*\circ \sigma$ and $g^*$ are equimeasurable (\cite[Chapter~2, Proposition~7.2]{BS})\textemdash and so are the functions $g$ and $g^*$\textemdash we have
\begin{equation*}
\mu_f(\lambda) = \mu(\{x\in\R\colon g^*(\sigma(x)) > \lambda \}) = |\{t>0\colon g^*(t) > \lambda\}| = \mu_g(\lambda)
\end{equation*}
for every $\lambda > 0$. Hence $f$ and $g$ are equimeasurable. Furthermore, we have
\begin{align*}
\int_{\R} f(x)^p v(x) \d{\mu(x)} &= \int_{\R} g^*(\sigma(x))^p v_*(\sigma(x)) \d{\mu(x)}\\
&= \int_0^{\mu(\R)} g^*(t)^p v_*(t) \d{t}.
\end{align*}
Here we used \eqref{thm:hull_of_weighted_Lp:v=v_down_star_circ_sigma_finite_measure} and the fact that $\sigma\colon\R \to (0, \mu(\R))$ is a measure-preserving transformation. Hence \eqref{thm:hull_of_weighted_Lp:almost_equal_norms} with $\varepsilon = 0$ follows.

Finally, if the condition {(i)} of Theorem~\ref{thm:our_Ryff} with $f$ replaced by $v$ is satisfied, we have $S=T$. Now, if $S=T=\infty$, we are in the same position as in the preceding paragraph; hence \eqref{thm:hull_of_weighted_Lp:almost_equal_norms} holds with $\varepsilon = 0$. If $S=T < \infty$, we defined $\R_2$ as
\begin{equation*}
\R_2 = \{x\in\R\colon v(x) = S\}
\end{equation*}
instead of \eqref{thm:hull_of_weighted_Lp:R2_def}, and proceed as in the case $S < \infty$ and $\kf[v](S) < \infty$. The key difference is that \eqref{thm:hull_of_weighted_Lp:lower_bound} is now an equality. Consequently, \eqref{thm:hull_of_weighted_Lp:almost_equal_norms} follows with $\varepsilon = 0$. If the condition {(i)} of Theorem~\ref{thm:our_Ryff} with $f$ replaced by $v$ is not satisfied but {(ii)} of Theorem~\ref{thm:our_Ryff} is, we have $S = T$ and $\kf[v](S) = \infty$. In this case, we already proved that \eqref{thm:hull_of_weighted_Lp:almost_equal_norms} holds with $\varepsilon = 0$.

It remains for us to prove (iii). Assume that $v_*\equiv0$. Note that this is equivalent to the fact that $\kf[v](s) = \infty$ for every $s\in(0,\infty)$. In particular, we necessarily have $\mu(\R) = \infty$. Let $g\in\Mpl(\R, \mu)$ be such that $g^*(t) < \infty$ for every $t\in(0,\infty)$, and let $\varepsilon > 0$. 
In addition, assume that $g^*(t) > 0$ for every $t\in(0,\infty)$. Since $\kf[v] = \infty$ on $(0,\infty)$, obviously we have
	$
		\kf[v](\varepsilon^p/(2^{j+1}g^*(j)^p)) = \kf[v](\varepsilon^p/(2g^*(2^{-j})^p)) = \infty
	$
for every $j\in\N$. Thanks to this and the fact that $\mu$ is nonatomic, it is possible to 
find, inductively, mutually disjoint \mbox{$\mu$-measurable} sets $E_j, F_j\subset \R$, $j\in\N$, satisfying the following conditions:
\begin{align*}
E_j &\subset \left\{ x\in\R\colon v(x) < \frac{\varepsilon^p}{2^{j+1}g^*(j)^p} \right\} \setminus \bigcup_{k = 1}^{j-1}(E_k \cup F_k), \\
\mu(E_j) &= 1,\\
F_j &\subset \left\{ x\in\R\colon v(x) < \frac{\varepsilon^p}{2g^*(2^{-j})^p} \right\} \mathbin{\Big\backslash} \left(\bigcup_{k = 1}^{j-1}(E_k \cup F_k) \cup E_j\right), \\
\intertext{and}
\mu(F_j) &= 2^{-j}.
\end{align*}
By \cite[Chapter~2, Proposition~7.4]{BS}, there are measure-preserving transformations $\sigma_j\colon E_j \to (j, j+1)$ and $\widetilde{\sigma}_j\colon F_j \to (2^{-j}, 2^{-j + 1})$. For each $j\in\N$, set
\begin{equation*}
f_j(x) = g^*(\sigma_j(x))\chi_{E_j}(x) + g^*(\widetilde{\sigma}_j(x))\chi_{F_j}(x),\ x\in\R.
\end{equation*}

Now, define the function $f\in\Mpl(\R, \mu)$ as
\begin{equation*}
f = \sum_{j = 1}^\infty f_j.
\end{equation*}
We claim that $f$ is equimeasurable with $g$ and satisfies \eqref{thm:hull_of_weighted_Lp:epsilon_norm}. As for the equimeasurability, we have, for every $\lambda > 0$,
\begin{align*}
\mu_f(\lambda) &= \sum_{j = 1}^\infty \mu(\{x\in E_j\colon g^*(\sigma_j(x)) > \lambda \}) + \sum_{j = 1}^\infty \mu(\{x\in F_j\colon g^*(\widetilde{\sigma}_j(x)) > \lambda \}) \\
&= \sum_{j = 1}^\infty |\{t\in (j, j+1)\colon g^*(t) > \lambda \}| + \sum_{j = 1}^\infty |\{t\in (2^{-j}, 2^{-j + 1})\colon g^*(t) > \lambda \}| \\
&= |\{t > 0\colon g^*(t) > \lambda \}| = \mu_g(\lambda).
\end{align*}
Here we used the fact that $E_j, F_j$ are mutually disjoint and $\sigma_j$, $\widetilde{\sigma}_j$ are measure-preserving. As for the validity of \eqref{thm:hull_of_weighted_Lp:epsilon_norm}, note that
\begin{align*}
\int_R f(x)^p v(x) \d{\mu}(x) &= \sum_{j = 1}^\infty \int_{E_j} g^*(\sigma_j(x))^p v(x) \d{\mu}(x) + \sum_{j = 1}^\infty \int_{F_j} g^*(\widetilde{\sigma}_j(x))^p v(x) \d{\mu}(x) \\
&\leq \sum_{j = 1}^\infty g^*(j)^p \frac{\varepsilon^p}{2^{j+1}g^*(j)^p} \mu(E_j) + \sum_{j = 1}^\infty g^*(2^{-j})^p \frac{\varepsilon^p}{2g^*(2^{-j})^p} \mu(F_j) \\
& = \frac{\varepsilon^p}{2} + \frac{\varepsilon^p}{2}.
\end{align*}
Therefore, $f$ has the desired properties. This finishes the proof of (iii) under the additional assumption that $g^*(t) > 0$ for every $t\in(0,\infty)$. It is easy to modify the preceding steps for the situation when there is $a\in(0, \infty)$ such that $g^*(t) = 0$ for every $t\geq a$, and so we omit it.
\end{proof}

\begin{remark}
Under the hypotheses of Theorem~\ref{embleblor}, if  $v_*\in B_p$ \cite{MR989570}, i.e., if 
$$
\int_r^\infty\frac{v_*(t)}{t^p}\d t\le C\frac1{r^p}\int_0^r v_*(t)\d t, \quad\text{for every } r>0,
$$
then  $\Lambda^p(\R, v_*)$ is an r.i.\  Banach function space, when $1<p<\infty$ \cite{S:90,MR2308059}, and hence 
\begin{equation}\label{bpbanach}
L^p(\R, v)\subset\Lambda^p(\R, v_*)\subset (L^1+L^\infty)(R,\mu).
\end{equation}
For example, in $(\mathbb R, \d x)$, if $v(x)=|x|^\alpha$, $0\le\alpha<p-1$, then $v_*(t)\approx t^\alpha\in B_p$, and $\Lambda^p(\mathbb R, t^\alpha)=L^{p/(\alpha+1),p}(\mathbb R)$, the classical Lorentz space, is the r.i.\  Banach function  hull of $L^p(\mathbb R,|x|^\alpha)$.
\end{remark}

As a consequence of  Theorem~\ref{embedding}, we can now conclude the following result:

\begin{corollary}
If $p \in (1 ,\infty)$ and $v\colon \R \to [0, \infty)$ is a $\mu$-measurable function such that $v_*\in B_p$, then  $v^{-1/(p-1)}\in 
(L^1 + L^\infty)(\R, \mu).$

\begin{proof}
As in \eqref{bpbanach}, since $\Lambda^p(\R, v_*)\subset (L^1+L^\infty)(R,\mu)$ and $\d\nu=v\, \d\mu$, from \eqref{charctnorm} we obtain that
\begin{align*}
 \sup_{\substack{E\subset \R \\ \mu(E) \le 1}}\bigg\|\frac1v\bigg\|_{L^{p'}(E, v\d\mu)}
 &=\sup_{0<t\le 1}\sup_{\mu(E)=t}\bigg(\int_Ev^{-1/(p-1)}(x)\d\mu(x)\bigg)^{1/{p'}}\\
 &=\sup_{0<t\le 1}\bigg(\int_0^t(v^{-1/(p-1)})^*(s)\d s\bigg)^{1/{p'}}=\|v^{-1/(p-1)}\|_{(L^1 + L^\infty)(\R, \mu)}^{1/{p'}}< \infty.
 \end{align*}\qedhere
 \end{proof}
\end{corollary}

\end{document}